\documentclass[a4paper,12pt]{article}
\usepackage{amsmath}
\usepackage{amsfonts}
\usepackage{latexsym}
 \usepackage[colorlinks=true]{hyperref}
\hypersetup{urlcolor=blue, citecolor=red}

\newtheorem{theorem}{Theorem}[section]
\newtheorem{corollary}{Corollary}[section]

\newtheorem{proposition}{Proposition}[section]

\newtheorem{definition}[theorem]{Definition}
\newtheorem{remark}{Remark}[section]

\newcommand{\la}[0]{{\langle}}
\newcommand{\ra}[0]{{\rangle}}
\newcommand{\beq}[0]{\begin{equation}}
\newcommand{\eeq}[0]{\end{equation}}
\newcommand{\uu}[0]{\boldsymbol{u}}
\newcommand{\un}[0]{\boldsymbol{u}_n}
\newcommand{\ww}[0]{\textbf{w}}
\newcommand{\VV}[0]{\textbf{V}}
\newcommand{\HH}[0]{\textbf{H}}
\newcommand{\vv}[0]{\boldsymbol{v}}
\newcommand{\tow}[0]{\rightharpoonup}
\newcommand{\RR}[0]{\mathbb{R}}
\newcommand{\R}[0]{\mathbb{R}}
\newcommand{\mint}[0]{\int_\Omega}
\newcommand{\tint}[0]{\int_0^t}
\newcommand{\Tint}[0]{\int_0^T}
\newcommand{\epsi}[0]{\varepsilon}
\newcommand{\ff}[0]{\varphi}

\newcommand{\EE}[0]{\mathcal{E}(\ff(t))}
\newcommand{\EEk}[0]{\mathcal{E}(\ff_k(t))}

\newcommand{\tmu}[0]{\tilde{\mu}}

\newcommand{\fn}[0]{\ff_n}

\newcommand{\mn}[0]{\mu_n}
\newcommand{\Ltv}[0]{L^2(0,T;V)}

\newcommand{\rd}[0]{\mathbb{R}^d}
\newcommand{\barf}[0]{\bar{\varphi}}

\newcommand{\dert}[0]{\frac{d}{dt}}
\newcommand{\nn}[0]{\boldsymbol{n}}
\newcommand{\CH}[0]{\ff_t + \nabla \cdot (\uu \ff) =\Delta \mu}
\newcommand{\locpot}[0]{\mu=-\Delta\ff+ F'(\ff)}
\newcommand{\nlchK}[0]{\ff_t + \nabla \cdot (\uu \ff) =\Delta \mu}
\newcommand{\Potential}[0]{\mu=a\ff-J\ast\ff+ F'(\ff)}

\numberwithin{equation}{section}

\title{On the nonlocal Cahn-Hilliard-Brinkman \\ and Cahn-Hilliard-Hele-Shaw systems}

\author{
{\sc Francesco Della Porta}\\
Mathematical Institute, University of Oxford\\
Oxford OX2 6GG, UK\\
\textit{francesco.dellaporta@maths.ox.ac.uk}
\\
\\
{\sc Maurizio Grasselli}\\
Dipartimento di Matematica, Politecnico di Milano\\
Milano 20133, Italy \\
\textit{maurizio.grasselli@polimi.it}
}

\begin{document}
\maketitle


%
%


\begin{abstract}
The phase separation of an isothermal incompressible binary fluid in a porous medium can be described by
the so-called Brinkman equation coupled with a convective Cahn-Hilliard (CH) equation.
The former governs the average fluid velocity $\mathbf{u}$, while the latter rules evolution of $\varphi$, the difference
of the (relative) concentrations of the two phases. The two equations are known as the Cahn-Hilliard-Brinkman (CHB) system.
In particular, the Brinkman equation is a Stokes-like equation with a forcing term (Korteweg force) which is proportional to $\mu\nabla\varphi$, where  $\mu$ is the chemical potential. When the viscosity vanishes, then the system becomes the Cahn-Hilliard-Hele-Shaw (CHHS) system. Both systems have been studied from the theoretical and the numerical viewpoints. However, theoretical results on the CHHS system are still rather incomplete. For instance, uniqueness of weak solutions is unknown even in 2D. Here we replace the usual CH equation with its physically more relevant nonlocal version. This choice allows us to prove more about the corresponding nonlocal CHHS system. More precisely, we first study well-posedness for the CHB system, endowed with no-slip and no-flux boundary conditions. Then, existence of a weak solution to the CHHS system is obtained as a limit of solutions to the CHB system. Stronger assumptions on the initial datum allow us to prove uniqueness for the CHHS
  system. Further regularity properties are obtained by assuming additional, though reasonable, assumptions on the interaction kernel. By exploiting these properties, we provide an estimate for the difference between the solution to the CHB system and the one to the CHHS system with respect to viscosity.
\end{abstract}



\begin{section}{Introduction}
\label{EUCHB}

The phenomenon of phase separation of incompressible binary fluids in a porous medium can be modeled by means of a diffuse interface approach. Consider a mixture of two fluids occupying a
bounded domain $\Omega\subset \mathbb{R}^d$, $d=2,3$, for any time $t\in (0,T)$, $T>0$, denote by $\varphi$ the difference of the fluid (relative) concentrations and by $\uu$  the (averaged) fluid velocity.
Assuming that the two fluids have the same constant density, the resulting model is the so-called Cahn-Hilliard-Brinkman (CHB) system (see, e.g., \cite{NYT, Schmuck2013})
\begin{equation}
\label{LCHB}
\begin{cases}
\CH \\
\locpot \\
-\nabla\cdot(\nu\nabla\uu)+\eta\uu+\nabla p=\mu\nabla\ff+\boldsymbol{h}  \\
\nabla \cdot \uu=0
\end{cases}
\end{equation}
in $\Omega\times(0,T)$, $T>0$.
Here $\nu>0$ is the viscosity coefficient, $\eta>0$ the fluid permeability and $p$ is the fluid pressure. Other constants
are supposed to be one for simplicity. The mobility is also assumed to be constant and equal to one, while
$F$ stands for a double well potential accounting for phase separation.
The average velocity $\uu$ obeys a modified Darcy's law proposed by H.C. Brinkman in 1947 (see \cite{Br}).

System \eqref{LCHB} endowed with no-slip and no-flux boundary conditions
has been analyzed from the numerical viewpoint in \cite{CSW} (see also \cite{DFW}). Some theoretical results can be found
in \cite{BCG}, where well-posedness in a weak setting as well as longtime
behavior of solutions (i.e., existence of the global attractor and convergence to a unique equilibrium)
have been investigated. Another interesting issue is the analysis of behavior of solutions when $\nu$ goes to zero.
Indeed when $\nu=0$ system \eqref{LCHB} becomes the so-called Cahn-Hilliard-Hele-Shaw (CHHS) model which
is used, for instance, to describe tumor growth dynamics (see, e.g., \cite{LTZ} and references therein, cf. also \cite{DFRSS}).
This model presents several technical difficulties (cf. \cite{LTZ,WW,WZ}, see also \cite{FW,DFW,Wise} for numerical schemes).
For instance, uniqueness of weak solutions is an open issue even in dimension two, as well as the existence of a global strong solution
in dimension three for sufficiently general initial data (see \cite{LTZ}). Existence
of a global weak solution to the CHHS system is obtained in \cite{BCG} as limit of solutions to system \eqref{LCHB}  (see also \cite[Thm.2.4]{FW} for an existence result). In the same paper, the difference of (strong) solutions to \eqref{LCHB} and the CHHS system is estimated with respect to $\nu$ and to the initial data in dimension two. Most of the quoted papers deal with a regular potential $F$, that is, $F$ is defined on the whole real line (however, see \cite{DFRSS} for a singular potential).

In this contribution we want to analyze a nonlocal variant of \eqref{LCHB} which is obtained by replacing the standard Cahn-Hilliard (CH) equation by its nonlocal version. More precisely, we consider the following nonlocal CHB system
\begin{equation}
\label{CHB}
\begin{cases}
\nlchK \\
\Potential \\
-\nabla\cdot(\nu(\ff)\nabla\uu)+\eta\uu+\nabla p=\mu\nabla\ff+\boldsymbol{h}  \\
\nabla \cdot \uu=0
\end{cases}
\end{equation}
in $\Omega\times(0,T)$. Here the viscosity may depend on $\ff$, while $J:\rd \to \mathbb{R}$ is a suitable interaction kernel and $a(x)=\int_\Omega J(x-y)dy$. This system is endowed with boundary and initial conditions
\begin{equation}
\label{IBC}
\begin{cases}
\displaystyle\frac{\partial\mu}{\partial \nn}=0 \qquad		&\text{on $\partial\Omega\times(0,T)$}\\
\uu=\textbf{0} \qquad		&\text{on $\partial\Omega\times(0,T)$}\\
\ff(0)=\ff_0 \qquad		&\text{in $\Omega$}.
\end{cases}
\end{equation}

We recall that the nonlocal CH equation can be justified in a more rigorous way from the physical viewpoint (cf. \cite{GL3}, see also \cite{GL1,GL2}). Also, the standard CH equation can be interpreted as an approximation of the nonlocal one. The nonlocal CH equation has been analyzed in a number of papers, under various assumptions on the potential $F$ and on the mobility (see, e.g., \cite{CH forte, dpg, MR, NLCHAnalysis1, NLCHExponential, NLCHAA1, NLCHAA2,RS}, cf. also \cite{GLWW,GWW} for the numerics). In addition, a series of papers have recently been devoted to the so-called Cahn-Hilliard-Navier-Stokes (CHNS) system in its nonlocal version (cf. \cite{grass,grass trattore,FG2,FGG,FGK,FGR,FRS}). Adapting the techniques devised in \cite{grass}, we can prove existence of a global weak solution to \eqref{CHB}--\eqref{IBC}. Its uniqueness (for constant viscosity) also holds in dimension three. However, the main goal is the analysis of the vanishing viscosity case where the limi
 t problem is
\beq
\label{CHHS}
\begin{cases}
\nlchK \qquad\\
\eta\uu+\nabla p=\mu\nabla\ff+\boldsymbol{h} \\
\nabla\cdot\uu=0
\end{cases}
\eeq
in $\Omega\times(0,T)$, i.e. the nonlocal CHHS system, subject to the boundary and initial conditions
\beq
\begin{cases}
\label{BCIC-HS}
\displaystyle\frac{\partial\mu}{\partial \nn}=0 \qquad		&\text{on $\partial\Omega\times(0,T)$}\\
\uu\cdot \nn = 0 \qquad		&\text{on $\partial\Omega\times(0,T)$}\\
\ff(0)=\ff_0 \qquad		&\text{in $\Omega$}.
\end{cases}
\eeq

As in \cite{BCG}, we can prove that a solution to \eqref{CHHS}--\eqref{BCIC-HS} can be obtained as a limit of solutions to \eqref{CHB}--\eqref{IBC}. In addition, uniqueness holds when $\ff_0$ is bounded (and so is $\ff$). Here we take advantage of the fact that the nonlocal CH equation is essentially a second-order equation and not a fourth-order equation like in the standard CHHS system. Then, further reasonable assumptions on $J$ allow us to establish some regularity properties of the solutions. These properties help us to estimate the difference, with respect to $\nu$ and the initial data, between a solution to \eqref{CHB}--\eqref{IBC} and a solution to the CHHS system.

The plan of this paper goes as follows. Notation, assumptions and statements of the main results are contained in Section \ref{s:preliminar}. Results concerning existence and regularity for \eqref{CHB}--\eqref{IBC} are proven in Section \ref{s:existenceCHB}. Existence of a weak solution to \eqref{CHHS}--\eqref{BCIC-HS} is demonstrated in Section \ref{s:existenceCHHS}. Section \ref{s:uniqueness} deals with uniqueness and continuous dependence on data for both problems.
The final Section \ref{s:comparison} is essentially devoted to obtain the estimate of the difference of the solutions to \eqref{CHB}--\eqref{IBC} and \eqref{CHHS}--\eqref{BCIC-HS}.

\section{Functional setup and main results}
\label{s:preliminar}
\subsection{Notation}
We set $H:=L^2(\Omega)$ and $V:=H^1(\Omega)$.
We denote by $\Arrowvert\,\cdot\,\Arrowvert$ and $(\cdot\,,\cdot)$ the norm and the scalar product in $H$, respectively, while
$\la\,\cdot\,\ra$ stands for the duality between $V^\prime$ and $V$.
For every $\ff\in V'$ we denote by $\bar \ff$ the average of $\ff$ over $\Omega$, namely $\bar \ff=|\Omega|^{-1}\la\ff,1\ra$.
Then we define
\[
V_2=\Bigl\{v \in H^2(\Omega):\, \frac{\partial v}{\partial\nn}=0\, \text{on}\, \partial\Omega\Bigr\}.
\]
The linear operator $A=-\Delta: V_2\subset H \to H$ with dense domain is self-adjoint and non-negative. Moreover, it is strictly positive on $V_0=\{\psi\in V\,:\, \bar{\psi} =0\}$ and it maps $V_0$ isomorphically into $V_0'=\{\psi\in V'\,:\,\la \psi,1\ra  =0\}$.
We will also set $$\|\psi - \bar{\psi}\|_{r}=\Arrowvert\, A^{r/2}(\psi - \bar{\psi})\,\Arrowvert $$ for every $r\in\R$.
Observe that the norm $\|\cdot\|_\#$ defined as
\[
\|x\|_\#:=\Bigl(\|x-\bar{x}\|^2_{-1}+\bar{x}^2 \Bigr)^\frac{1}{2},
\]
is equivalent to the usual norm of $V'$.

Besides, let $\mathcal{V}$ be the space of divergence-free test functions defined by
\begin{equation*}
    \mathcal{V} = \{ \vv \in C^{\infty}_{0}(\Omega, \mathbb{R}^{d}) \, : \, \nabla \cdot \vv = 0 \}.
\end{equation*}
We shall use the following canonical spaces (see, e.g., \cite[Chapter~I]{Temam2001})
$$
    {\boldsymbol H}    = \overline{\mathcal{V}}^{H^{d}} \quad\text{ and }\quad
    {\boldsymbol V}    = \{ \vv \in V^d \, : \, \nabla \cdot \vv = 0 \}.
$$
Recall that  $\vv \in {\boldsymbol V}$ yields $\vv |_{\partial \Omega} = \mathbf{0}$, while $\vv \in {\boldsymbol H}$ is such that $\vv \cdot \nn = 0$ on $\partial \Omega$. We will still use $(\cdot\,,\cdot)$ and $\la\,\cdot\,\ra$ to denote the scalar product in $\HH$ and the duality between ${\boldsymbol V}^\prime$ and ${\boldsymbol V}$, respectively.

Finally, $c$ will indicate a generic nonnegative constant depending on $\Omega$, $J,\, F,$ and $\boldsymbol{h}$ at most. Instead,
$N$ will stand for a generic positive constant which has further dependence on $T$ and/or on some norm of $\ff_0$.
The value of $c$ and $N$ may vary even within the same line.

\subsection{Assumptions}
Following~\cite{CH forte} and~\cite{grass} (cf. also~\cite{BCG})
we introduce the following assumptions.
\begin{enumerate}
	\item[(H0)]
		  $\Omega\subset\mathbb{R}^d$, $d=2,3$, is open, bounded and connected with a smooth boundary.
	\item[(H1)]
		$J\in W^{1,1}(\rd)$ satisfies
		\[
         J(x)=J(-x), \qquad a(x):=\mint J(x-y)\,dy \geq 0, \;\text{a.e.}\, x\in\Omega.
		\]
	\item[(H2)]
		$F\in C^{2,1}_{loc}(\R)$ and there exists $c_0 > 0$ such that
		\[
F''(s)+a(x)\geq c_0, \quad \forall s\in \R,\; \text{a.e.}\, x \in\Omega.
\]
	\item[(H3)]
		There exist $c_1>0$, $c_2>0$ and $q>0$ if $d=2$, $q\geq\frac{1}{2}$ if $d=3$ such that
		$$
		F''(s)+a(x)\geq c_{9}|s|^{2q}-c_{10},\qquad\forall s\in\R,\text{ a.e. }x\in\Omega.
		$$
	\item[(H4)]
		There exist $c_3>0$ and $p\in (1,2]$ such that
		\[
		|F'(s)|^p\leq c_4(|F(s)|+1),\quad \forall s\in \R.
		\]
	\item[(H5)]
		$\eta\in L^\infty(\Omega)$ and
$$\eta(x)\geq 0, \quad \text{a.e. } x\in\Omega.$$
	\item[(H6)]
		$\nu$ is locally Lipschitz on $\RR$ and there exist $\nu_0,\nu_1>0$ such that
		$$\nu_0\leq\nu(s)\leq\nu_1, \quad \forall{s\in\RR}.$$		
	\item[(H7)]
		$\boldsymbol{h}\in L^2(0, T; \VV')$.
\end{enumerate}

\noindent

\begin{remark}
\label{perturbation}
Assumption (H2) implies that the potential $F$ is a quadratic
perturbation of a strictly convex function. Indeed $F$ can be
represented as
\beq
\label{F coerciva}
F(s)=G(s)-\frac{a^*} {2}\,s^2
\eeq
with $G\in C^{2,1}(\R)$ strictly convex, since $G''\geq c_0$ in
$\Omega$. Here $a^*= \| a\| _{L^{\infty}(\Omega)}$ and observe that $a\in L^\infty(\Omega)$ derives from (H1).
\end{remark}	

\begin{remark}
\label{growth}
Since $F$ is bounded from below, it is easy to see that (H4) implies that $F$ has polynomial growth of order $p'$, where $p'\in[2,\infty)$ is the conjugate index to p. Namely there exist $c_4>0$ and $c_5\geq 0$ such that
\[
|F(s)|\leq c_4|s|^{p'}+c_5,\qquad \forall s\in\R.
\]
Besides, it can be shown that (H3) implies the existence of $c_6,c_7>0$ such that
\[
F(s)\geq c_6|s|^{2+2q}-c_5,\qquad \forall s\in\R.
\]
\end{remark}	

\begin{remark}
The usual double well potential  $F(s)=\frac{1}{4}(s^2-1)^2$ satisfies all the hypotheses on $F$.
\end{remark}	

\begin{remark}
\label{h4}
One easily realizes that (H4) implies
		\[
		|F'(s)|\leq c(|F(s)|+1),\quad \forall s\in \R;
		\]
furthermore (H3) implies that
		\[
		|F(s)|\leq F(s)+ 2\max\{0,c_2\},\quad \forall s\in \R.\\\\
		\]
\end{remark}	

\begin{remark}
Note that (H5) allows, in particular, $\eta=0$. Thus the so-called Cahn-Hilliard-Stokes system is also included
(see \cite{Schmuck2013}).
\end{remark}

\begin{remark}
The convective nonlocal CH equation can formally be rewritten as follows
\begin{equation*}
\ff_t = \nabla\cdot\bigl( (F''(\ff)+a)\nabla\ff\bigr)+\nabla\cdot\bigl(\nabla a\ff-\uu\ff\bigr)-\nabla J\ast \ff
\end{equation*}
from which the crucial role of (H2) is evident, namely, we are dealing with a convection-diffusion integrodifferential equation.
\end{remark}

\subsection{Statement of the main results}
Let us introduce the definition of weak solution to \eqref{CHB}--\eqref{IBC}.

\begin{definition}
\label{soluzione debole}
Let $T>0$ be given and let $\ff_0\in H$ be such that $F(\ff_0)\in L^1(\Omega)$.
A pair $(\ff,\uu)$ is a weak solution to \eqref{CHB}--\eqref{IBC} on $[0,T]$ if
\begin{align*}
&\ff\in C([0,T];H)\cap L^2(0,T;V)\\
&\ff_t \in L^2(0,T;V')\\
&\Potential \in \Ltv\\
&\uu \in L^2(0,T;\VV)
\end{align*}
and it satisfies
\begin{align}
\label{weak ch}
&\la  \ff_t,\,\psi\ra +(\nabla\mu,\,\nabla\psi)=(\uu\ff,\,\nabla\psi)
, \quad\forall \psi\in V,\quad\text{a.e. in }(0,T),\\
\label{weak B}
&(\nu(\ff)\nabla\uu,\,\nabla\vv)+(\eta\uu,\,\vv) = (\mu\nabla\ff,\,\vv)+\la\boldsymbol{h},\,\vv\ra, \quad\forall \vv\in \VV,\quad\text{a.e. in } (0,T),\\
&\label{ci}
\ff(0)=\ff_0, \quad\text{a.e. in } \Omega.
\end{align}
\end{definition}

\begin{remark}
Observe that if we choose $\psi=1$ in \eqref{weak ch} we obtain
\[
\dert \barf = 0.
\]
Thus the total mass of any weak solution is conserved.
\end{remark}

Global existence of a weak solution is given by
\begin{theorem}
\label{buona posizione}
Let $\ff_0\in H$ be such that $F(\ff_0)\in L^1(\Omega)$ and suppose that (H0)-(H7) are satisfied.
Then there exists a weak solution $(\ff,\uu)$ to \eqref{CHB}--\eqref{IBC}.
Furthermore, $F(\ff)$ is in $L^\infty(0,T;L^1(\Omega))$ and setting
\beq
\label{nlEnergy}
\mathcal{E}(\ff(t))=\frac{1}{4}\mint\mint J
(x-y)(\ff(x,t)-\ff(y,t))^2 \,dx\,dy + \mint F(\ff(x,t))\,dx.
\eeq
the following energy equality holds for almost every $t\in(0,T)$
\beq
\label{energy equality}
\dert \mathcal{E}(\ff(t)) + \|\nabla\mu\|^2 + \|\sqrt{\nu(\ff)}\nabla\uu\|^2+\|\sqrt\eta\uu\|^2 =\la  \boldsymbol{h} ,\uu\ra .
\eeq
\end{theorem}

Furthermore, we have
\begin{corollary}
\label{u in linf}
Let (H0)-(H6) hold. If $\boldsymbol{h}\in L^\infty(0,T;\VV')$ for some $T>0$. Then,
any weak solution $(\ff,\uu)$ to \eqref{CHB}--\eqref{IBC} is such that
 $$\ff\in L^4(0,T;L^4(\Omega)),\qquad \uu\in L^\infty(0,T;V).$$
\end{corollary}

Weak solutions can be regular provided $\ff_0$ is bounded. Indeed we have
\begin{proposition}
\label{ff in linf}
Let the assumptions of Theorem \ref{buona posizione} hold. If $\ff_0\in L^\infty(\Omega)$ then, any solution $(\ff,\uu)$ to problem \eqref{CHB} on $[0,T]$ corresponding to $\ff_0$ satisfies
$$\ff,\mu\in L^\infty(\Omega\times(0,T)).$$
In particular, we have
$$\|\ff\|_{L^\infty(\Omega\times(0,T))}\leq M,\qquad \|\mu\|_{L^\infty(\Omega\times(0,T))}\leq M,$$
for some $M>0$, independent of $\nu$ and $T$.
\end{proposition}

If the viscosity $\nu$ is constant then we have a continuous dependence estimate
\begin{proposition}
\label{stability}
Let hypotheses (H0)-(H5) hold. Suppose that $\nu$ is a positive constant and $\boldsymbol{h}\in L^\infty(0,T;\VV^\prime)$. Consider two weak solutions to \eqref{CHB}--\eqref{IBC}, namely $(\ff_1,\uu_1)$ and $(\ff_2,\uu_2)$, corresponding to the initial data $\ff_{1,0}$ and $\ff_{2,0}$, respectively. Here $\ff_{i,0}\in L^2(\Omega)$ and $F(\ff_{i,0})\in L^1(\Omega)$, $i=1,2$. Then there exists $N=N(T)>0$ such that, for any $t\in [0,T],$
\begin{align}
\label{chb dip continua}
\|\ff_1(t)-\ff_2(t)\|^2_\#+\tint\|\uu_1(y)-\uu_2(y)\|^2_\VV \,dy\leq N\bigl(\|\ff_{1,0}-\ff_{2,0}\|^2_\#+|\barf_{1,0}-\barf_{2,0}|\bigr).
\end{align}
In particular, \eqref{CHB}--\eqref{IBC} has a unique weak solution.
\end{proposition}

\medskip
\noindent
\textbf{The limit $\nu\to 0$.}
\label{EUCHHS}
As a second step in our analysis we study the limit of \eqref{CHB}--\eqref{IBC} with constant viscosity $\nu$,
as $\nu$ tends to 0. We recall that the resulting limit system is \eqref{CHHS}--\eqref{BCIC-HS} whose weak formulation is given by the following definition.

\begin{definition}
\label{soluzione debole CHHS}
Let $T>0$ be given and let $\ff_0\in L^\infty(\Omega)$.

A pair $(\ff,\uu)$ is a weak solution to \eqref{CHHS}--\eqref{BCIC-HS} on $(0,T)$ if
	\begin{align*}
	&\ff\in L^\infty(\Omega \times (0,T))\cap L^2(0,T;V)\\
	&\ff_t\in L^2(0,T;V')\\
	&\Potential \in \Ltv\\
	&\uu\in L^2(0,T;\HH)
	\end{align*}
and it satisfies
\begin{align}
\label{weak CHHS}
&\la \ff_t,\psi\ra +(\nabla\mu,\nabla\psi)=(u\ff,\nabla\psi), \quad\forall \psi\in V,\quad\text{a.e. in }(0,T),\\
\label{weak HS}
&(\eta\uu,\,\vv) = (\mu\nabla\ff,\,\vv)+(\boldsymbol{h},\,\vv), \quad\forall \vv\in \HH,\quad\text{a.e. in } (0,T),\\
\nonumber
&\ff(0)=\ff_0,\quad\text{a.e. in } \Omega.
\end{align}
\end{definition}

\medskip
\noindent
To analyze \eqref{CHHS}--\eqref{BCIC-HS} we replace assumption (H5) with the stronger
\begin{enumerate}
\item[(H8)] $\eta\in L^\infty(\Omega)$ and there exists $\eta_0>0$ such that
		$$\eta(x)\geq \eta_0,\quad \text{a.e.}\, x\in\Omega .$$
\end{enumerate}
Furthermore, for the sake of simplicity, we let $\boldsymbol{h}=\textbf{0}$.
Then we have the following existence theorem

\begin{theorem}
\label{buona posizione CHHS}
Let (H0)-(H4), (H8) hold and let $\ff_0\in L^\infty(\Omega)$.
Then, for any given $T>0$, if $\{\nu_k\}$ is a sequence of positive constants converging to $0$, the weak solution to \eqref{CHB}--\eqref{IBC} with $\nu=\nu_k$ converges, up to a subsequence, to a weak solution $(\ff,\uu)$ to \eqref{CHHS}--\eqref{BCIC-HS}. More precisely, we have
\begin{align*}
&\ff_k\to \ff \quad \text{strongly  in }  L^2(0,T;H) \;\\
&\uu_k\tow \uu \quad \text{weakly  in }  L^2(0, T; \HH)
\end{align*}
Furthermore, the following energy equality holds for almost any $t\in (0,T)$:
\beq
\label{energy equality CHHS}
\dert \mathcal{E}(\ff(t)) + \|\nabla\mu\|^2 + \|\sqrt\eta\uu\|^2 =0,
\eeq
where $\mathcal{E}$ is defined by \eqref{nlEnergy}.
\end{theorem}

Next corollary is related to further regularity in the case where $\eta$ is constant.
\begin{corollary}
\label{reg-HS}
Let the assumptions of Theorem \ref{buona posizione CHHS} hold and $\eta$ be a positive constant, then $\uu \in L^\infty(0,T;[L^p(\Omega)]^d)$ for each $p\geq 1$.
\end{corollary}

This fact allows us to prove uniqueness of the (weak) solution to \eqref{CHHS}--\eqref{BCIC-HS} for constant parameter $\eta$.
More precisely, we have
\begin{proposition}
\label{buona posizione CHHS 2}
Let the assumptions of Corollary \ref{reg-HS} hold. Consider two weak solutions to \eqref{CHHS}--\eqref{BCIC-HS}, namely $(\ff_1,\uu_1)$, $(\ff_2,\uu_2)$ corresponding to bounded initial data $\ff_{1,0}$, $\ff_{2,0}$, respectively. Then there exists $N=N(T)>0$ such that, for every $t\in [0,T]$,
$$
\|\ff_1-\ff_2\|^2_\#+\tint\|\uu_1-\uu_2\|^2_\HH
\leq N\bigl(\|\ff_{1,0}-\ff_{2,0}\|^2_\#+|\barf_{1,0}-\barf_{2,0}|\bigr).
$$
In particular, there exists a unique bounded weak solution to \eqref{CHHS}--\eqref{BCIC-HS}.
\end{proposition}

In case $J$ is more regular, we gain regularity also for the velocity field $\uu$. For the sake of completeness, we
first recall the definition of admissible kernel (see \cite[Definition 1]{BRB}).
\begin{definition}
A kernel $J\in W^{1,1}_{loc}(\R^d)$, $d=2,3$ is admissible if the following conditions are satisfied:
\begin{itemize}
\item $J\in C^3(\R^d \setminus \{0\})$;
\item $J$ is radially symmetric, i.e., $J(x) = \tilde J(|x|)$ with $\tilde J$ non-increasing;
\item $\tilde J''(r)$ and $\tilde J'(r)/r$ are monotone on $(0, r_0)$ for some $r_0 > 0$;
\item $|D^3J(x)|\leq C|x|^{-(d+1)}$ for some $C > 0$.
\end{itemize}
\end{definition}
Then we state the following regularity result
\begin{proposition}
\label{regu h1 u}
Let the assumptions of Theorem \ref{buona posizione CHHS} hold, $\eta$ be constant and $J$ be admissible or $J\in W^{2,1}$. Then
$$
\uu\in L^2(0,T;\VV).
$$
\end{proposition}

Thanks to the above regularity result we can obtain an estimate of the difference between a solution to \eqref{CHB}--\eqref{IBC} and a solution to \eqref{CHHS}--\eqref{BCIC-HS}. Indeed we have
\begin{theorem}
\label{t:close}
Let (H0), (H2)-(H4), (H8) hold. Suppose $\nu$, $\eta$ constant, $\mathbf{h}=\mathbf{0}$, and
$J$ either be admissible or $J\in W^{2,1}(\R^2)$. Take $\ff_0^\nu,\ff_0\in L^\infty(\Omega)$ and
$$R:=\sup_{\nu> 0}\{\|\ff_0^\nu\|_{L^\infty},\|\ff_0\|_{L^\infty}\}<\infty.$$
Let $(\ff_\nu,\uu_\nu)$ be the unique weak solution to \eqref{CHB}--\eqref{IBC} with initial datum $\ff_0^\nu$,
and $(\ff, \uu)$ the unique solution to \eqref{CHHS}--\eqref{BCIC-HS} with initial datum $\ff_0$.
Then, for any given $T>0$, there exists $C_{R,T}>0$ such that
$$\|\ff_\nu(t)-\ff(t)\|_\#^2+\int_0^t\|\uu_\nu(y)-\uu(y)\|^2\,d y\leq
\big(\|\ff_0^\nu-\ff_0\|_\#^2+|\barf_0^\nu-\barf_0|\big)e^{C_{R,T}}+C_{R,T}\nu,$$
for each $t\in [0,T]$. In particular, if $\ff_0^\nu=\ff_0$, then
$\ff_\nu\to \ff$ in $L^\infty(0,T;V')$ and in $L^2(0,T;H)$ as $\nu\to 0.$
\end{theorem}

\section{Existence and regularity for the CHB system}
\label{s:existenceCHB}
The first part of this section is devoted to prove Theorem \ref{buona posizione}.
Then, in the second part, the proofs of Corollary \ref{u in linf} and Proposition \ref{ff in linf} are given.
\subsection*{Proof of Theorem \ref{buona posizione}}
\noindent
The proof will be carried out by means of a Faedo--Galerkin approximation scheme, following closely~\cite{grass}. We first prove existence of a solution when $\ff_0\in V_2$ and $\boldsymbol{h}\in C([0,T]; \HH)$; then, by a density argument, we will recover the same result for any initial datum  $\ff_0\in H$ with $F(\ff_0)\in L^1(\Omega)$ and any $\boldsymbol{h}\in L^2([0,T]; \VV')$.

We consider the families $\{\psi_j\}_{j\in\mathbb{N}}\subset V_2$ and $\{\vv_j\}_{j\in\mathbb{N}}\subset \VV$ respectively
eigenvectors of $A+I: V_2 \to H$ and of the Stokes operator, which are both self-adjoint, positive and linear.
Let us define the $n$-dimensional subspaces
$\Psi_n:=\la \psi_1,...,\psi_n\ra $ and $\mathcal{W}_n:=\la \ww_1,...,\ww_n\ra $
with the related
orthogonal projectors on this subspace $P_n:=P_{\Psi_n}$ and $\tilde{P}_n:=P_{\mathcal{W}_n}$.
We then look for three functions of the following form:
\[
\ff_n(t)=\sum_{k=1}^{n} b_k^{(n)}(t)\psi_k, \qquad
\mu_n(t)=\sum_{k=1}^{n} c_k^{(n)}(t)\psi_k, \qquad
\uu_n(t)=\sum_{k=1}^{n} d_k^{(n)}(t)\ww_k
\]

\noindent
that solve the following discretized problem
\begin{align}
\label{eq approssimata}
&(\ff_n',\psi)+(\nabla\rho_n,\nabla\psi)
=(\un\ff_n,\nabla\psi)+(\nabla
J\ast\ff_n,\nabla\psi) \\
\label{eq approssimata 2}
&(\nu(\fn)\nabla\un,\,\nabla\ww)+(\eta\un,\,\ww) + (\fn\nabla\mn,\,\ww) = \la  \boldsymbol{h},\,\ww\ra \\
\label{defrho}
&\rho_n:=a(\cdot)\ff_n+F'(\ff_n),\\
\label{apprmu}
&\mu_n=P_n(\rho_n-J\ast\ff_n),\\
\label{ci approssimata}
&\ff_n(0)=\ff_{0n},
\end{align}
\noindent
for every $\psi\in\Psi_n$, every $\ww\in\mathcal{W}_n$ and where $\ff_{0n}:=P_n\ff_0$.

By using the definition of $\ff_n$, $\mn$ and $\un$, problem  \eqref{eq approssimata}--\eqref{ci approssimata} becomes equivalent to a Cauchy problem for a system of ordinary
differential equations in the $n$ unknowns $b_i^{(n)}$.
Thanks to (H2), the Cauchy-Lipschitz theorem yields that there exists a unique solution $b^{(n)}\in C^1([0,T^*_n];\R^n)$
for some maximal time $T^*_n\in(0,+\infty]$.

Let us show that $T^*_n=+\infty$, for all $n\geq 1$. Indeed, using $\psi=\mu_n$ as test function in \eqref{eq approssimata} and $\ww=\un$ in~\eqref{eq approssimata 2}
we get the following identity:
\beq
\label{es chb 0}
(\ff_n',\mu_n)+(\nabla\rho_n,\nabla\mu_n)+\|\sqrt{\nu(\fn)}\nabla\un\|^2+\|\sqrt{\eta}\un\|^2
=(\nabla
J\ast\ff_n,\nabla\mu_n) + \la  {\boldsymbol{h}},\un\ra .
\eeq
Let us first notice that
\begin{align}
\label{es chb 1}
&(\ff_n',\mu_n)=\frac{d}{dt}\Bigl( \frac{1}{4}\mint\mint
J(x-y)(\ff_n(x)-\ff_n(y))^2+\mint F(\ff_n) \Bigr),\\
\label{es chb 2}
&(\nabla \mu_n,\nabla P_n(J\ast\ff_n)) \leq
\frac{1}{4}\|\nabla\mu_n\|^2+\|\ff_n\|^2\|J\|_{W^{1,1}}^2,\\
\label{es chb 3}
&(\nabla J\ast\ff_n,\nabla\mu_n)\leq
\frac{1}{4}\|\nabla\mu_n\|^2+\|\ff_n\|^2\|J\|_{W^{1,1}}^2.
\end{align}

\noindent
\\
By means of (H3), we can deduce the existence of a positive constant $\alpha$ such that
\begin{align}
\nonumber
&\frac{1}{2}\mint\mint J(x-y) (\ff_n(x)-\ff_n(y))^2\,dx\,dy
+2\mint F(\ff_n)
\\
\label{es chb 4}
&=\|a\ff_n\|^2 + 2\mint F(\ff_n) - (\ff_n, J\ast \ff_n)
\geq \alpha\Bigl(
 \| \ff_n\|^2
+ \mint F(\fn)\Bigr)
- c.
\end{align}

\noindent
By using (H6) and Poincar\'e's inequality, it is easy to show that there exists $\beta>0$ such that
\beq
\label{es chb 5}
\beta \|\un\|^2_\VV\leq\|\sqrt{\nu}\nabla\un\|^2,
\eeq
and, on account of (H7,) we have
\beq
\label{es chb 6}
\la  \boldsymbol{h}, \un\ra  \leq c\|\boldsymbol{h}\|_{\VV'}^2+\frac{\beta}{2}\|\un\|_\VV^2.
\eeq

Let us now exploit~\eqref{es chb 1} in~\eqref{es chb 0} and integrate it with respect to
time between $0$ and $t\in(0,T_n^*)$. Taking \eqref{es chb 2}--\eqref{es chb 6} into account, we find
\begin{align}
\label{maxinter}
&\alpha\Bigl( \| \ff_n\|^2 + \mint F(\fn)\Bigr)
+\tint\Bigl(\frac{\beta}{2}\|\un(\tau)\|_\VV^2+\|\sqrt{\eta}\un(\tau)\|^2+\|\nabla\mu_n(\tau)\|^2\Bigr)\, d\tau \nonumber\\
&\leq M + K\tint \Bigl( \| \ff_n(\tau)\|^2 + \mint F(\fn(\tau))\Bigr) \, d\tau,
\end{align}
\noindent
which holds for all $t\in [0,T_n^*),$ where
\[
M=c\Bigl(1 + \|\ff_{0}\|^2 + \mint F(\ff_0)+ \|\boldsymbol{h} \| ^2_{L^2 (0,T;\VV') }\Bigr),
\]
and $K=2\|J\|^2_{W^{1,1}}.$
Here, we have used the fact that that $\ff_0$ and $\ff_{0,n}$ are supposed to
belong to $V_2$. We point out that $M$ and $K$ do not depend on $n$.

Thus, inequality  \eqref{maxinter} entails that $T_n^*=+\infty$, for all $n\geq 1$.
As a consequence, \eqref{eq approssimata}--\eqref{ci approssimata}
has a unique global-in-time solution. Furthermore, we obtain the
following estimates, holding for any given $0<T<+\infty$:
\beq
\label{fn limitata}
\|\ff_n\|_{L^\infty(0,T;H)}\leq N
\eeq
\beq
\label{mn limitata}
\|\nabla\mu_n\|_{L^2(0,T;H)}\leq N
\eeq
\beq
\label{Ffn limitata}
\|F(\ff_n)\|_{L^\infty(0,T;L^1(\Omega))}\leq N
\eeq
\beq
\label{un limitata}
\|\un\|_{L^2(0,T;\VV)}\leq \frac{N}{\sqrt{\nu_0}}
\eeq
where $N$ is independent of $n$. Observe that, in light of (H3), \eqref{Ffn limitata} implies
\beq
\label{fn limitata 2+2q}
\|\ff_n\|_{L^\infty(0,T;L^{2+2q}(\Omega))}\leq N
\eeq
Thanks to (H2), recalling  \eqref{apprmu}, we get
\begin{align*}
\frac{c_0}{4}\|\nabla\fn\|^2+\frac{1}{c_0}\|\nabla\mn\|^2&\geq
(a\nabla\ff_n+\ff_n\nabla a+F''(\ff_n)\nabla\ff_n-\nabla J*\ff_n,\nabla\ff_n)\\
&\geq c_0\|\nabla
\fn\|^2-2\|\nabla J\|_{L^1}\|\nabla
\fn\|\|\ff_n\|\\
&\geq \frac{c_0}{2}\|\nabla
\fn\|^2-c\|\fn\|^2,
\end{align*}
thus \eqref{fn limitata} and \eqref{mn limitata} yield
\beq
\label{fn in V}
\|\fn\|_{L^2(0,T;V)}\leq N.
\eeq

\noindent
The next step is to deduce a (uniform) bound for $\mn$ in $\Ltv$.
Thanks to Remark~\ref{h4} and to the identity
$$
(P_n(-J\ast\fn+a\fn),1)=(-J\ast\fn+a\fn,1)=0
$$
we get
\beq
\label{mn minore di F}
\biggl | \mint\mn \biggr|
= \bigl|(F'(\fn),1)\bigr|
\leq \mint \bigl|F'(\fn)\bigr|
\leq c\mint F(\fn)+c\leq N.
\eeq
The Poincar\'e inequality implies
\beq
\biggl\|\mn-\frac{1}{|\Omega|}\mint\mn \biggl\|  \leq  c\|\nabla\mn\|,
\eeq
and from \eqref{mn limitata} and \eqref{mn minore di F} we deduce that
\beq
\label{mn in V}
\|\mn\|_{\Ltv}\leq N.
\eeq

\noindent
Observe now that, calling $\tilde\rho_n=P_n\rho_n$,
\[
\|\tilde\rho_n\|^2_V=\|\mn+P_n(J\ast\fn)\|_V^2
\leq 2\|\mn\|_V^2+ 2(\|J\|^2_{L^{1}}+\|\nabla
J\|^2_{L^{1}})\|\fn\|^2,
\]
so that from \eqref{mn in V} we immediately get
\beq
\label{rho in V}
\|\tilde\rho_n\|_{\Ltv}\leq N.
\eeq
Furthermore, recalling  \eqref{defrho} and invoking (H4), we obtain
\[
\|\rho_n\|_{L^p }
\leq c a^*\|\fn\|+\|F'(\fn)\|_{L^p }
\leq  cN+c\Bigl( \mint |F(\fn)| \Bigr)^{1/p}\leq N,
\]
which yields the bound
\beq
\label{rho in Lp}
\|\rho_n\|_{L^\infty(0,T; L^p(\Omega))}\leq N.
\eeq
We finally provide an estimate for the sequence $\fn'$. We take a generic test function $\psi\in V$ and we
write it as $\psi=\psi_1+\psi_2$, where $\psi_1=P_n\psi\in\Psi_n$ and $\psi_2=\psi-\psi_1\in \Psi_n^\perp$.
It is easy to see that
\beq
\label{stima termine rho}
|(\nabla\rho_n,\nabla\psi_1)|
\leq \|\nabla\tilde\rho_n\| \|\nabla\psi_1\|
\leq \|\nabla\tilde\rho_n\| \|\nabla\psi\|_V,\\
\eeq
and
\begin{align}
\label{es chb 11}
&|(\un\fn,\nabla\psi_1)|
\leq \|\un\|_{[L^{\frac{2+2q}{q}}]^d} \|\nabla\psi_1\| \|\fn\|_{L^{2+2q}}
\leq  N \|\un\|_\VV \|\psi\|_V,
&\quad d=2,\\
\label{es chb 11 bis}
&|(\un\fn,\nabla\psi_1)|
\leq \|\un\|_{[L^6]^d} \|\nabla\psi_1\| \|\fn\|_{L^3}
\leq  N \|\un\|_\VV \|\psi\|_V,
&\quad d=3.
\end{align}
\noindent
By using Young's lemma we infer
\beq
\label{stima termine J}
\Bigl| \mint\nabla J\ast\fn\nabla\psi_1 \Bigr|
\leq \|\psi\|_V \|\nabla J\|_{L^1}\|\fn\|
\leq  N \|\nabla J\|_{L^1}\|\psi\|_V.
\eeq
\noindent
From \eqref{eq approssimata}, owing to \eqref{stima termine rho}-\eqref{stima termine J}, we have that
\beq
\label{stima per ff' 1}
|(\fn',\psi)|
\leq N(1+\|\nabla\rho_n\|+\|\un\|_\VV)\|\psi\|_V,
\eeq
which gives
\beq
\label{fn' in V'}
\|\fn'\|_{L^2(0,T;V')}\leq N,
\eeq
owing to \eqref{un limitata} and \eqref{rho in V}.
\noindent
Collecting estimates  \eqref{fn limitata},  \eqref{fn in V},  \eqref{mn in V}--\eqref{rho in Lp}, \eqref{fn' in V'}, we find
\begin{align*}
&\ff\in L^\infty(0,T;L^{2+2q}(\Omega))\cap\Ltv\cap H^1(0,T;V'),\\
&\mu\in\Ltv,\\
\nonumber
&\tilde{\rho}\in\Ltv,\\
\nonumber
&\rho\in L^\infty (0,T; L^p(\Omega)),\\
\nonumber
&\uu\in L^2 (0,T; \VV),
\end{align*}
such that, up to a subsequence,
\begin{align}
\label{conv 0}
&\fn\tow \ff \quad \text{weakly*  in } L^\infty(0,T;H),\\
\label{fn conv}
&\fn\tow \ff \quad \text{weakly  in } \Ltv,\\
\label{strong conv}
&\fn\to \ff \quad \text{strongly  in }  L^\gamma(0,T;H) \;\text{and a.e.  in } \Omega\times(0,T),\\
\label{mn conv}
&\mn\tow \mu \quad \text{weakly  in } \Ltv,\\
\label{rhon conv}
&\tilde\rho_n\tow \tilde\rho \quad \text{weakly  in }  \Ltv,\\
\label{rho conv}
&\rho_n \tow \rho \quad \text{weakly*  in } L^\infty (0,T; L^p(\Omega)),\\
\label{fn' conv}
&\fn'\tow \ff_t \quad \text{weakly  in }  L^2(0, T; V'),\\
\label{un conv}
&\un\tow \uu \quad \text{weakly  in }  L^2(0, T; \VV).
\end{align}
\noindent
Here $\gamma = 2+2q$ if $d=2$, $\gamma = \min\{{2+2q},{4}\}$ if $d=3$. We now pass to the limit in \eqref{eq approssimata}--\eqref{ci approssimata} in order to prove that
$(\ff,\uu)$ is a weak solution to CHB system according to Definition \ref{soluzione debole}.
First of all, from the pointwise convergence \eqref{strong conv} we have $\rho_n\to a\ff+F'(\ff)$ almost everywhere in $\Omega\times (0,T)$, therefore from \eqref{rho conv} we have $\rho=a\ff+F'(\ff)$.
Now, for every $\phi\in\Psi_j$, every $j\leq n$ with $j$ fixed and for every $\chi\in C_0^\infty(0,T)$, we have that
\[
\Tint(\rho_n ,\phi)\chi(t)=\Tint(\tilde\rho_n,\phi)\chi(t).
\]
Passing to the limit in this equation, using \eqref{rhon conv} and \eqref{rho conv}, and on account of the density of $\{\Psi_j\}_{j\geq 1}$ in $H$, we get $\tilde{\rho}(\cdot,\ff)=\rho(\cdot,\ff) = a\ff+F'(\ff)$.
Moreover, since $\mn=P_n(\rho_n-J\ast\fn)$, then, for every $\phi\in\Psi_j$, every $k\leq j$ with $j$ fixed and for every $\chi\in C_0^\infty(0,T)$, there holds
\[
\Tint(\mu_n(t),\phi)\chi(t)dt=\Tint(\rho_n-J\ast\fn,\phi)\chi(t)dt.
\]
By passing to the limit in the above identity, and using the convergences \eqref{strong conv}, \eqref{mn conv} and \eqref{rho conv}, we eventually get
\[
\Potential=\rho-J\ast\ff.
\]
It still remains to pass to the limit in \eqref{eq approssimata} and \eqref{eq approssimata 2} in order to recover \eqref{weak ch}, \eqref{weak B} and initial condition \eqref{ci}. This can be obtained in a standard way, so we refer the reader to ~\cite[Proof of Theorem 1]{grass} where all the technicalities are detailed. In order to conclude to proof,
let us now assume that $\ff_0\in H$ with $F(\ff_0)\in L^1(\Omega)$ and $\mathbf{h}\in L^2(0,T; \VV')$. In this case, we first choose an approximating sequence of initial data
$\ff_{0n}\in V_2$ such that $\ff_{0n}\to \ff_0$ in $H$, and a sequence $\mathbf{h}_n\in C(0,T; \HH)$ in such a way that $\mathbf{h}_n\to \mathbf{h}$ in $L^2(0,T; \VV')$. Then, arguing as in~\cite[Proof of Theorem 1]{grass} the existence of a solution to \eqref{CHB}--\eqref{IBC} is obtained by passing to the limit $n\to\infty$.  In particular, on account of  \eqref{es chb 4}-\eqref{es chb 6}, we find that $F(\ff)\in L^\infty(0,T;L^1(\Omega))$.

We are left to prove the energy identity \eqref{energy equality}. Let us take $\psi=\mu(t)$ in equation \eqref{weak ch}. This yields
\beq
\label{ee 1}
\la  \ff_t, \mu\ra  +\|\sqrt\nu\nabla\uu\|^2+\|\sqrt\eta\uu\|^2+\|\nabla\mu\|^2=\la  \boldsymbol{h},\uu\ra .
\eeq
By arguing as in~\cite[proof of Corollary 2]{grass}, one can prove the identity
\[
\la  \ff_t, \mu\ra  = \la  \ff_t, a\ff+F'(\ff)-J\ast\ff\ra  =\dert \EE
\]
which holds for almost every $t>0$. Thus \eqref{energy equality} follows directly from \eqref{ee 1}.
\hfill$\Box$

\subsection*{Proof of Corollary \ref{u in linf}}
We recall that a standard application of the Gagliardo-Nirenberg inequality gives
\begin{align*}
&\|\ff\|_{L^4}\leq \|\ff\|^{1/2}\|\nabla\ff\|^{1/2},\qquad & d=2,\\
&\|\ff\|_{L^4}\leq \|\ff\|_{L^3}^{1/2}\|\nabla\ff\|^{1/2},\qquad & d=3.
\end{align*}
On account of Theorem \ref{buona posizione}, we have $\ff\in L^2(0,T;V)$. Moreover, owing to (H3), we have
$\ff\in L^\infty(0,T;L^2(\Omega))$ if $d=2$ and $\ff\in L^\infty(0,T;L^{3}(\Omega))$ if $d=3$. Then we easily deduce
\[
\int_0^T\|\ff\|^4_{L^4}\leq N\int_0^T\|\ff\|^2_V\leq N.
\]
In order to prove the estimate for $\uu$, let us first recall the following identity (see \cite[Proof of Thm.~2]{FGG})
\beq
\label{urka}
(\mu \nabla\ff, \uu)=(\nabla J*\ff, \ff\uu)-(\frac12\nabla a\ff^2, \uu).
\eeq
Thanks to \eqref{urka}, equation \eqref{weak B} with $\vv=\uu$ can be rewritten as follows
\beq
\label{reg 1}
\|\sqrt{\nu}\nabla\uu\|^2+\|\sqrt{\eta}\uu\|^2=(\nabla J\ast\ff,\,\ff\uu)-\frac12(\nabla a\ff^2,\,\uu)+\la  \boldsymbol{h},\,\uu\ra.
\eeq
Observe now that
\begin{align*}
&(\nabla J\ast\ff,\,\ff\uu)-\frac12(\nabla a\ff^2,\,\uu)\\
\leq& \left(\frac12\|\nabla a\|_{L^\infty}+\|\nabla J\|_{L^1}\right)\|\ff\|\|\ff\|_{L^{2+2q}}\|\uu\|_{L^{\frac{2+2q}{q}}}
,\quad &d=2,\\
&(\nabla J\ast\ff,\,\ff\uu)-\frac12(\nabla a\ff^2,\,\uu)\\
\leq& \left(\frac12\|\nabla a\|_{L^\infty}+\|\nabla J\|_{L^1}\right)\|\ff\|\|\ff\|_{L^3}\|\uu\|_{L^6}
,\quad &d=3.
\end{align*}
and, as $\ff\in L^\infty(0,T;L^{2+2q}(\Omega))$ when $d=2$, we obtain
$$
(\nabla J\ast\ff,\,\ff\uu)-\frac12(\nabla a\ff^2,\,\uu)\leq N \|\uu\|_\VV.
$$
On the other hand we get (cf.~(H6))
$$
\|\sqrt{\nu}\nabla\uu\|^2\geq \nu_0\|\nabla\uu\|^2\geq c\|\uu\|_\VV^2.
$$
Hence, by (H8) and \eqref{reg 1}, we end up with
\[
c\|\uu\|_\VV^2\leq N \|\uu\|_\VV
\]
which yields $\uu\in L^\infty(0,T;\VV)$.
\hfill$\Box$


\subsection*{Proof of Proposition \ref{ff in linf}}
In order to prove that $\ff\in L^\infty(\Omega\times(0,T))$ we can use a Moser-Alikakos type argument (see \cite[Proof of Thm.~3]{FGR} for the details). The boundedness of $\mu$ follows from its definition by comparison.
\hfill$\Box$



\begin{section}{Existence and regularity for CHHS system}
\label{s:existenceCHHS}

\subsection*{Proof of Theorem~\ref{buona posizione CHHS}}

Let $(\ff_k,\uu_k)$ be the solution of problem~\eqref{CHB} with $\nu=\nu_k$, thus satisfying~\eqref{nlEnergy}. Therefore, for every $k\geq 1$ we have
\[
\EEk + \tint \bigl( \|\nabla\mu_k\|^2+\|\sqrt\nu\nabla\uu\|^2+\|\sqrt\eta\uu\|^2\bigr)=\mathcal{E}(\ff_0)
\]
and thanks to~\eqref{es chb 4} it is possible to deduce~\eqref{fn limitata}--\eqref{Ffn limitata} and
\begin{align}
\label{un limitata VV hs}
&\|\uu_k\|_{L^2(0,T;\VV)}\leq \frac{N}{\sqrt\nu_k}\\
\label{un limitata HH hs}
&\|\uu_k\|_{L^2(0,T;\HH)}\leq N.
\end{align}
Furthermore, by arguing as in the proof of Theorem~\ref{buona posizione}, it is possible to recover~\eqref{fn in V} and~\eqref{mn in V}. Then from Proposition \ref{ff in linf} we deduce the following bound
\beq
\label{fn limitata in linf hs}
\|\ff_k\|_{L^\infty(\Omega\times(0,T))}\leq N.
\eeq
Also, we observe that
\beq
\label {ftk in l2 1 hs}
(\nabla\mu_k,\,\nabla\psi)\leq \|\nabla\mu\|\|\nabla\psi\|
\eeq
and (see~\eqref{fn limitata in linf hs})
\beq
\label {ftk in l2 2 hs}
(\uu_k\ff_k,\nabla\psi)\leq\|\uu_k\|_\HH\|\ff\|_{L^\infty}\|\nabla\psi\|.
\eeq
By exploiting~\eqref{ftk in l2 1 hs}--\eqref{ftk in l2 2 hs} in~\eqref{weak ch} we deduce~\eqref{fn' in V'} by comparison. We recall that $N$ does not depend neither on $k$ nor on $\nu_k$. Summing up, we deduce the existence of
\begin{align*}
&\ff\in L^\infty(\Omega\times(0,T))\cap\Ltv\cap H^1(0,T;V'),\\
&\mu\in\Ltv,\\
&\uu\in L^2 (0,T; \HH),
\end{align*}
such that, up to a subsequence,
\begin{align}
\label{hs conv 0}
&\ff_k\tow \ff \quad \text{weakly*  in } L^\infty(\Omega\times(0,T)),\\
\label{hs fn conv}
&\ff_k\tow \ff \quad \text{weakly  in } \Ltv,\\
\label{hs strong conv}
&\ff_k\to \ff \quad \text{strongly  in }  L^2(0,T;L^{\beta}(\Omega)) \;\text{and a.e.  in } \Omega\times(0,T),\\
\label{hs mn conv}
&\mu_k\tow \mu \quad \text{weakly  in } \Ltv,\\
\label{hs fn' conv}
&\ff_k'\tow \ff_t \quad \text{weakly  in }  L^2(0, T; V'),\\
\label{hs un conv}
&\uu_k\tow \uu \quad \text{weakly  in }  L^2(0, T; \HH).
\end{align}
Here $\beta$ is such that $\frac{1}{2}=\frac{1}{d+\epsi}+\frac{1}{\beta}$ for some $\epsi>0$.\\

It is now possible to pass to the limit as $k\to\infty$ in the weak formulation of \eqref{CHB}--\eqref{IBC}. We will do that restricting ourselves to the case $\psi\in W^{1,d+\epsi}(\Omega)\subset V$ in \eqref{weak ch} and then recovering the fact that \eqref{weak CHHS} holds for every $\psi\in V$ by a density argument.
Some attention is needed when passing to the limit in the viscous term of the Brinkman equation; as a matter of fact we have
\[
\nu_k (\nabla\uu_k,\,\nabla\vv)\leq\nu_k \|\nabla\uu_k\|\|\nabla\vv\|\leq \sqrt{\nu_k} N\|\nabla\vv\|
\]
which tends to $0$ as $\nu_k\to 0$. The convective term can be treated as follows:
\[
\int_t^{t+r} (\uu_k\ff_k-\uu\ff,\,\nabla\psi)=\int_t^{t+r} (\uu_k(\ff_k-\ff),\,\nabla\psi)+\int_t^{t+r} ((\uu_k-\uu)\ff,\,\nabla\psi)
\]
where $r\geq 0$ is arbitrary. Here the second term vanishes thanks to the boundedness of $\ff$ and~\eqref{hs un conv}. The first one goes to 0 thanks to~\eqref{un limitata HH hs},~\eqref{hs strong conv} and the fact that
\[
\int_t^{t+r} (\uu_k(\ff_k-\ff),\,\nabla\psi)\leq \|\ff-\ff_k\|_{ L^2(0,T;L^\beta)}\|\uu_k\|_{L^2(0, T;\HH)}\|\nabla\psi\|_{L^{d+\epsi}}.
\]
Finally, we can pass to the limit into the the Korteweg force since, for every $r\geq 0$, we have
\[
\int_t^{t+r} (\nabla\mu_k\ff_k-\nabla\mu\ff,\,\vv)=\int_t^{t+r} (\nabla\mu_k (\ff_k-\ff),\,\vv)+\int_t^{t+r} (\nabla(\mu_k-\mu)\ff,\,\vv)
\]
and the second term goes to $0$ thanks to the boundedness of $\ff$ and~\eqref{hs mn conv}, while the first one vanishes thanks to~\eqref{mn in V} and~\eqref{hs strong conv} and the inequality
\[
\int_t^{t+r} (\nabla\mu_k (\ff_k-\ff),\,\vv)\leq \|\vv\|_\VV\|\mu_k\|_{\Ltv}\|\ff-\ff_k\|_ {L^2(0,T;L^3)}.
\]
It is easy to see that~\eqref{weak HS} makes sense also for every $\vv\in\HH$. Furthermore, thanks to~\eqref{ftk in l2 2 hs} we can deduce that \eqref{weak CHHS} holds also for every $\psi\in V$ by a density argument. Thus, we showed that there is a subsequence of $(\ff_k,\uu_k)$ converging to a $(\ff,\uu)$ which is a weak solution to \eqref{CHHS}--\eqref{BCIC-HS}.
\hfill$\Box$

\subsection{Proof of Corollary \ref{reg-HS}}
As shown in \eqref{urka}, we can rewrite \eqref{weak HS} as
\beq
\label{eq con kort new}
\eta(\uu,\vv)=(\nabla J*\ff, \ff\vv)-(\frac12\nabla a\ff^2, \vv),\qquad \text{ a.e. in }[0,T],\forall\vv\in \HH.
\eeq
On account of Lemma 2.1 in \cite{LTZ} we can deduce
$$
\|\uu\|_{[L^p]^d}\leq c(\|(\nabla J*\ff)\ff\|_{L^p}+\|\nabla a\ff^2\|_{L^p}).
$$
Furthermore, from Theorem \ref{buona posizione CHHS} we have $\ff\in L^\infty(0,T;\Omega)$, which, thanks to (H1), leads to  $\uu \in L^\infty(0,T; L^p(\Omega))$ for each $p\geq 1$.

\subsection{Proof of Proposition \ref{regu h1 u}}
As $\eta$ is constant we can take advantage of Lemma 2.1 in \cite{LTZ} and, rewriting the Korteweg force as in \eqref{eq con kort new} we can write
\beq
\label{est -1 h1 u}
\|\uu\|_\VV\leq c (\|\ff \nabla J*\ff\|_V+\frac12\|\nabla a\ff^2\|_V).
\eeq
As $\ff\in L^\infty(0,T;\Omega)$, from (H1) we can easily deduce that
\beq
\label{est 0 h1 u}
\|\ff \nabla J*\ff\|\leq c \|\ff\|^2_{L^\infty},\qquad\|\nabla a\ff^2\|\leq c\|\ff\|^2_{L^\infty}.
\eeq
Besides, we have
\beq
\label{est 1 h1 u}
\|\nabla(\ff \nabla J*\ff)\| \leq \| (\nabla J*\ff )\otimes \nabla\ff \|+\|\ff \nabla^2 J*\ff\|\leq c(\|\ff\|_{L^\infty}\|\nabla\ff\|+\|\ff\|^2_{L^\infty})
\eeq
and
\beq
\label{est 2 h1 u}
\|\nabla(\nabla a\ff^2)\| \leq \|\nabla^2 a\ff^2\|+2\|\ff \nabla a \otimes \nabla \ff\|\leq c(\|\ff\|_{L^\infty}\|\nabla\ff\|+\|\ff\|^2_{L^\infty}).
\eeq
Therefore, collecting \eqref{est -1 h1 u}-\eqref{est 2 h1 u} we finally conclude the proof of the proposition.

\end{section}

\section{Continuous dependence and uniqueness}
\label{s:uniqueness}
\subsection*{Proof of Proposition \ref{stability}}
\label{unique sect}
Let $(\ff_1,\textbf{u}_1)$ and $(\ff_2,\textbf{u}_2)$ be two weak solutions to the system \eqref{CHB}--\eqref{IBC}
corresponding to  $\ff_{1,0}$ and $\ff_{2,0}$, respectively.
Here $N>0$ will denote a generic constant depending on $T$ and $\|\ff_{i,0}\|$, $i=1,2$.\\
Setting $\ff=\ff_1-\ff_2$, $\tmu=\mu(\ff_1)-\mu(\ff_2)$ and $\uu=\uu_1-\uu_2$, we have
\begin{align}
\label{uniq chb2}
&\la \ff_{t},\psi\ra  + (\nabla\tmu,\nabla\psi) = (\textbf{u}\ff_1,\nabla \psi)+(\textbf{u}_2\ff,\nabla \psi),
\quad\forall \psi\in V,\quad\text{a.e. in }(0,T),\\
\label{uniq chb2 bis}
&\nu(\nabla\uu,\nabla\vv)+(\eta\uu,\vv)=(\tmu\nabla\ff_1,\vv)+(\mu_2\nabla\ff,\vv),
\quad\forall \vv\in \VV,\quad\text{a.e. in } (0,T),\\
\label{ICdiff}
&\ff(0)=\ff_{1,0}-\ff_{2,0}, \quad\text{a.e. in } \Omega.
\end{align}
Choosing $\psi=1$ we readily obtain that $\barf(t)=\ff(0)$ for all $t\in [0,T]$. On account of this, let us take $\psi=(-\Delta)^{-1}(\ff-\barf)$  in
\eqref{uniq chb2} and find
\begin{align}
\label{uniq chb3}
&\frac{1}{2}\dert \|\ff-\barf\|^2_{-1}+ (\tmu,\ff-\barf) = I_1+I_2,
\end{align}
where
$$I_1=(\textbf{u}\ff_1,\nabla (-\Delta)^{-1}(\ff-\barf)),\qquad I_2=(\textbf{u}_2\ff,\nabla (-\Delta)^{-1}(\ff-\barf)).$$

\noindent
Furthermore, taking  $\vv=\uu$ in~\eqref{uniq chb2 bis}, we get
$$
\nu\|\nabla\uu\|^2+\|\sqrt{\eta}\uu\|^2=(\tmu\nabla\ff_1,\uu)+(\mu_2\nabla\ff,\uu).
$$
After standard computations in light of \eqref{urka}, we obtain
$$
(\tmu\nabla\ff_1,\uu)+(\mu_2\nabla\ff,\uu)=(\nabla J\ast\ff_1,\,\ff\uu)+(\nabla J\ast\ff,\,\ff_2\uu)-\frac12(\nabla a (\ff_1+\ff_2),\,\ff\uu).
$$
If $d=2$, since $\ff_i\in L^\infty(0,T;L^{2+2q}(\Omega))$, $i=1,2$, then we obtain
\begin{align}
\label{uniq 6}
\nonumber
&(\tmu\nabla\ff_1,\uu)+(\mu_2\nabla\ff,\uu)\\
\nonumber
&\leq\max\left\{\frac12\|\nabla a\|_{L^\infty},\,\|\nabla J\|_{L^1})\right\} \|\uu\|_{[L^{\frac{2+2q}{q}}]^d}(\|\ff_1\|_{L^{2+2q}}+\|\ff_2\|_{L^{2+2q}})\|\ff\|\\
&\leq N\|\uu\|_{\VV}\|\ff\|.
\end{align}
Analogously, if $d=3$, recalling that $\ff_i\in L^\infty(0,T;L^3(\Omega))$, $i=1,2$, we deduce
\begin{align*}
(\tmu\nabla\ff_1,\uu)+(\mu_2\nabla\ff,\uu)&\leq\max\left\{\|\nabla a\|_{L^\infty},\,\|\nabla J\|_{L^1}\right\} \|\uu\|_{[L^6]^d}(\|\ff_1\|_{L^3}+\|\ff_2\|_{L^3})\|\ff\|\\
&\leq N\|\uu\|_{\VV}\|\ff\|.
\end{align*}
Observe now that
$$
\nu\|\nabla\uu\|_\VV^2+\|\sqrt{\eta}\uu\|^2\geq c\|\uu\|^2_\VV
$$
gives
\beq
\label{uniq 7}
\|\uu\|_\VV\leq N\|\ff\|.
\eeq
Let us now estimate the terms in the differential equality \eqref{uniq chb3}.
In order to estimate $(\tmu,\ff-\barf)$ we argue as in \cite[proof of Proposition 2.1]{dpg} to deduce
\beq
\label{uniq 1}
(a\ff+F'(\ff_1)-F'(\ff_2),\ff-\barf)\geq \frac{7c_0}{8}\|\ff\|^2-c\barf^2-N|\barf|
\eeq
and
\beq
\label{uniq 2}
(J\ast\ff,\,\ff-\barf)\leq \frac{c_0}{8}\|\ff\|^2 + c\|\ff-\barf\|^2_\#.
\eeq
On the other hand, we have
\begin{align*}
&(\textbf{u}\ff_1,\nabla (-\Delta)^{-1}(\ff-\barf))\leq \|\ff_1\|_{L^{2+2q}}\|\uu\|_{[L^{\frac{2+2q}{q}}]^d}\|\ff-\barf\|_\#,&\qquad d=2,\\
&(\textbf{u}\ff_1,\nabla (-\Delta)^{-1}(\ff-\barf))\leq \|\ff_1\|_{L^3}\|\uu\|_{[L^6]^d}\|\ff-\barf\|_\#,&\qquad d=3.
\end{align*}
implying
\beq
\label{uniq 3}
I_1\leq N\|\uu\|_{\VV}\|\ff-\barf\|_\#.
\eeq
Concerning $I_2$, suppose $d=2$ first and observe that
\[
(\textbf{u}_2\ff,\nabla (-\Delta)^{-1}(\ff-\barf))\leq\frac{c_0}{16}\|\ff\|^2+c\|\uu_2\|^2_{[L^4]^d}\| \nabla(-\Delta)^{-1}(\ff-\barf)\|^2_{L^4}
\]
and
\begin{align*}
\| \nabla(-\Delta)^{-1}(\ff-\barf)\|^2_{L^4}&\leq c\|\nabla (-\Delta)^{-1}(\ff-\barf)\|\| \nabla(-\Delta)^{-1}(\ff-\barf)\|_V\\
&\leq c\|\ff-\barf\|\|\ff-\barf\|_\#.
\end{align*}
Thus, on account of Corollary \ref{u in linf}, we get
\begin{align*}
\label{uniq 4}
(\textbf{u}_2\ff,\nabla (-\Delta)^{-1}(\ff-\barf))&\leq\frac{c_0}{8}\|\ff\|^2+\|\uu_2\|^4_{\VV}\|\ff-\barf\|^2_{\#}+c\barf^2
\\\nonumber&\leq\frac{c_0}{8}\|\ff\|^2+N\|\ff-\barf\|^2_{\#},
\end{align*}
so that
\beq
\label{uniq 4}
I_2\leq \frac{c_0}{8}\|\ff\|^2+N\|\ff-\barf\|^2_{\#}.
\eeq
Inequality~\eqref{uniq 4} can also be proved in the case $d=3$ by considering
\[
(\textbf{u}_2\ff,\nabla (-\Delta)^{-1}(\ff-\barf))\leq\frac{c_0}{16}\|\ff\|^2+\|\uu_2\|^2_{[L^6]^d}\|\nabla (-\Delta)^{-1}(\ff-\barf)\|^2_{L^3},
\]
and observing that
\begin{align*}
\|\nabla (-\Delta)^{-1}(\ff-\barf)\|^2_{L^3}&\leq c\|\nabla (-\Delta)^{-1}(\ff-\barf)\|\|\nabla (-\Delta)^{-1}(\ff-\barf)\|_V\\
& \leq c\|\ff-\barf\|\|\ff-\barf\|_\#.
\end{align*}
Collecting \eqref{uniq 1}--\eqref{uniq 4}, we deduce from  \eqref{uniq chb3} the differential inequality
\beq
\label{uniq 9a}
\frac{1}{2}\dert \|\ff-\barf\|^2_{-1}+ \frac{c_0}{4}\|\ff\|^2\leq
N\|\uu\|_{\VV}\|\ff-\barf\|_\#+
N\|\ff-\barf\|_\#^2
+c\barf^2-N|\barf|.
\eeq
Taking \eqref{uniq 7} into account, we deduce
\beq
\label{uniq 9}
\frac{1}{2}\dert \|\ff-\barf\|^2_{\#}+\frac{c_0}{8} \|\ff\|^2\leq N\|\ff-\barf\|^2_{\#}+N|\barf|
\eeq
and Gronwall's lemma yields
$$\|\ff_1(t)-\ff_2(t)\|^2_\#\leq N\bigl(\|\ff_{1,0}-\ff_{2,0}\|^2_\#+|\barf_{1,0}-\barf_{2,0}|\bigr).
$$
The estimate for $\uu$ follows from \eqref{uniq 7} by integrating \eqref{uniq 9} on $[0,t]$, $t\in (0,T]$.
\hfill$\Box$

\subsection*{Proof of Proposition \ref{buona posizione CHHS 2}}
We argue in the same way as in the Proof of Proposition~\ref{stability}.
However, in this case we take advantage of the inequality
\beq
\label{hs uniq 2}
(\eta\uu,\,\uu)\geq\eta_0\|\uu\|^2.
\eeq
Moreover, we observe that~\eqref{uniq 6} can be replaced by
\begin{align}
\label{uniq_6bis}
\nonumber
(\tmu\nabla\ff_1,\,\uu)+&(\mu_2\nabla\ff,\,\uu)=
(\nabla J\ast\ff_1,\,\ff\uu)+(\nabla J\ast\ff,\,\ff_2\uu)-(\nabla a (\ff_1+\ff_2),\,\ff\uu)\\
&\leq \max\{\|\nabla a\|_{L^\infty},\,\|\nabla J\|_{L^1}\}\|\uu\|\bigl(\|\ff_1\|_{L^\infty}+\|\ff_2\|_{L^\infty}\bigr)\|\ff\|.
\end{align}
Leveraging on the fact that $\ff_1$ and $\ff_2$ are bounded, we obtain
\[
\|\uu\|\leq N\|\ff\|.
\]
Consider now \eqref{uniq chb3}. Instead of controlling $I_1$ as in~\eqref{uniq 3}, we obtain
\beq
\label{hs uniq 1}
I_1=(\uu\ff_1,\,\nabla(-\Delta)^{-1}(\ff-\barf))\leq N \|\uu\| \|\ff-\barf\|_\#
\eeq
Also, exploiting the estimates for $\uu$ and arguing as in the proof of Proposition \ref{stability}, thanks to Corollary \ref{reg-HS} we have
\begin{align*}
I_2&=(\textbf{u}_2\ff,\nabla (-\Delta)^{-1}(\ff-\bar \ff))\\
&\leq\frac{c_0}{8}\|\ff\|^2+\|\uu_2\|^4_{L^{2d}}\|\ff-\barf\|^2_{\#}+c\barf^2
\leq\frac{c_0}{8}\|\ff\|^2+N\|\ff-\barf\|^2_{\#}.
\end{align*}
Thus we can still prove inequality \eqref{uniq 9a} and the proof can be completed arguing as above.
\end{section}
\hfill$\Box$

\section{Convergence of solutions as $\nu\to 0$}
\label{s:comparison}
In this section we prove Theorem \ref{t:close}.

\subsection{Proof of Theorem \ref{t:close}}
We first define $\psi=\ff_\nu-\ff$, $\tilde\mu=\mu(\ff_\nu)-\mu(\ff)$ and $\vv=\uu_\nu-\uu$.
Let us now take $\vv$ in the weak formulation of the equation for $\vv$. Adding
$-\nu ( \nabla \uu, \nabla \vv ) $ to both sides of the resulting identity, we get
\[
\nu \| \nabla \vv \|^2 + \|\sqrt{\eta}\vv\|^2=( \tilde \mu\nabla \ff_\nu ,\vv) +(  \mu \nabla \psi,\vv)
-\nu(  \nabla \uu,\nabla \vv)
\]
Since
$$-\nu(  \nabla \uu,\nabla \vv) \leq \nu\|\nabla\uu\|^2+\nu\|\nabla\vv\|^2$$
we obtain
$$
\eta\|\vv\|^2\leq |(  \tilde \mu\nabla \ff_\nu ,\vv) +(   \mu \nabla \psi,\vv) |
+\nu\|\nabla\uu\|^2.
$$
Reasoning as in \eqref{uniq_6bis} we find
\begin{align*}
|(\tmu\nabla\ff_\nu,\,\vv)+(\mu\nabla\psi,\,\vv)|
&\leq \max{(\|\nabla a\|_{L^\infty},\,\|\nabla J\|_{L^1})}\|\vv\|\bigl(\|\ff_\nu\|_{L^\infty}+\|\ff\|_{L^\infty}\bigr)\|\psi\|\\
&\leq C\|\vv\|\|\psi\|,
\end{align*}
hence
$$
\eta\|\vv\|^2\leq C\|\vv\|\|\psi\|+\nu\|\nabla\uu\|^2.
$$
Note that this implies
\begin{equation}
\label{lavV}
\|\vv\|\leq \frac{C}{\eta}\|\psi\|+\frac{\sqrt{\nu}}{\sqrt{\eta}}\|\nabla\uu\|.
\end{equation}
On the other hand, we have
\begin{equation*}
\frac{1}{2}\dert \|\psi-\bar\psi\|^2_{-1}+ (\tmu,\psi-\bar\psi) = I_1+I_2,
\end{equation*}
where
$$I_1=(\textbf{v}\ff_\nu,\nabla (-\Delta)^{-1}(\psi-\bar\psi)),\qquad I_2=(\textbf{u}\psi,\nabla (-\Delta)^{-1}(\psi-\bar\psi)).$$
Now, by arguing as in proof of Proposition \ref{buona posizione CHHS 2} and exploiting boundedness of $\uu$ we deduce

\begin{align*}
&\frac{1}{2}\dert \|\psi-\bar{\psi}\|^2_{-1}+ \frac{c_0}{4}\|\psi\|^2\leq
N\|\vv\|\|\psi-\bar{\psi}\|_\#+
N\|\psi-\bar{\psi}\|_\#^2
+c\bar{\psi}^2+N|\bar{\psi}|.
\end{align*}
Thus, taking \eqref{lavV} into account, we end up with
\begin{equation}
\label{grou}
\frac{1}{2}\dert \|\psi-\bar{\psi}\|^2_\#+ \frac{c_0}{8}\|\psi\|^2\leq
N\|\psi-\bar{\psi}\|_\#^2
+N|\bar{\psi}|+N\nu\|\nabla\uu\|^2.
\end{equation}
An application of the Gronwall lemma on $[0,T]$, on account of Proposition \ref{regu h1 u} provides
$$\|\ff_\nu(t)-\ff(t)\|_\#^2\leq
\big(\|\ff_0^\nu-\ff_0\|_\#^2+|\barf_0^\nu-\barf_0|\big)e^{C_T}+C_T\nu.
$$
Now a further integration of \eqref{grou}, and \eqref{lavV} complete the proof.
\hfill$\Box$

\section*{Acknowledgments} The work of the first author was supported by the Engineering and Physical Sciences Research Council [EP/L015811/1]. The second author is member of the Gruppo Nazionale per l'Analisi Matematica, la Probabilit\`{a} e le loro Applicazioni (GNAMPA) and of the Istituto Nazionale di Alta Matematica (INdAM).


\end{document}